\documentclass[10pt,twoside]{article}
\usepackage{psfig,amsfonts,amssymb}
\newcommand{\Dt}{\Delta}
\newcommand{\DtL}{\Delta_2^{(L)}}

\newcommand{\La}{\Lambda}
\newcommand{\ZZ}{\mathbb Z}
\newcommand{\FF}{\mathbb F}
\newcommand{\QQ}{\mathbb Q}
\newcommand{\CC}{\mathbb C}

\newcommand{\RR}{\mathbb R}

\def\firstpage{1}\def\lastpage{1000}

\makeatletter
\def\magnification{\afterassignment\m@g\count@}
\def\m@g{\mag=\count@\hsize6.5truein\vsize8.9truein\dimen\footins8truein}
\makeatother

\font\caps=cmcsc10                    
\font\Caps=cccsc10 scaled \magstep1   
\font\scaps=cmcsc8

\makeatletter
\@specialpagefalse\headheight=8.5pt
\def\DocMath{{\def\th{\thinspace}\scaps Doc.\th Math.\th J.\th DMV}}
\renewcommand{\@oddfoot}{\hfill\scaps Documenta Mathematica $\cdot$ Extra
        Volume ICM  1998 $\cdot$ \number\firstpage--\lastpage\hfill}
\renewcommand{\@evenfoot}{\ifnum\thepage>\lastpage\hfill\scaps
    Documenta Mathematica $\cdot$ Extra Volume ICM  1998 $\cdot$ \hfill\else\@oddfoot\fi}%
\renewcommand{\@evenhead}{%
    \ifnum\thepage>\lastpage\rlap{\thepage}\hfill%
    \else\rlap{\thepage}\slshape\leftmark\hfill\caps\SAuthor\hfill\fi}%
\renewcommand{\@oddhead}{%
    \ifnum\thepage=\firstpage{\DocMath\hfill\llap{\thepage}}%
    \else{\slshape\rightmark}\hfill\caps\STitle\hfill\llap{\thepage}\fi}%
\makeatother

\def\TSkip{\bigskip}
\newbox\TheTitle{\obeylines\gdef\GetTitle #1
\ShortTitle  #2
\SubTitle    #3
\Author      #4
\ShortAuthor #5
\EndTitle
{\setbox\TheTitle=\vbox{\baselineskip=20pt\let\par=\cr\obeylines%
\halign{\centerline{\Caps##}\cr\noalign{\medskip}\cr#1\cr}}%
	\copy\TheTitle\TSkip\TSkip%
\def\next{#2}\ifx\next\empty\gdef\STitle{#1}\else\gdef\STitle{#2}\fi%
\def\next{#3}\ifx\next\empty%
    \else\setbox\TheTitle=\vbox{\baselineskip=20pt\let\par=\cr\obeylines%
    \halign{\centerline{\caps##} #3\cr}}\copy\TheTitle\TSkip\TSkip\fi%
\centerline{\caps #4}\TSkip\TSkip%
\def\next{#5}\ifx\next\empty\gdef\SAuthor{#4}\else\gdef\SAuthor{#5}\fi%
\catcode'015=5}}\def\Title{\obeylines\GetTitle}
\def\Abstract{\begingroup\narrower
    \parskip=\medskipamount\parindent=0pt{\caps Abstract. }}
\def\EndAbstract{\par\endgroup\TSkip}

\long\def\MSC#1\EndMSC{\def\arg{#1}\ifx\arg\empty\relax\else
     {\par\narrower\noindent%
     1991 Mathematics Subject Classification: #1\par}\fi}

\long\def\KEY#1\EndKEY{\def\arg{#1}\ifx\arg\empty\relax\else
	{\par\narrower\noindent Keywords and Phrases: #1\par}\fi\TSkip}

\newbox\TheAdd\def\Addresses{\vfill\copy\TheAdd\vfill
    \ifodd\number\lastpage\vfill\eject\phantom{.}\vfill\eject\fi}
{\obeylines\gdef\GetAddress #1
\Address #2
\Address #3
\Address #4
\EndAddress
{\def\xs{6truecm}
\setbox0=\vtop{{\obeylines\hsize=\xs#1\par}}\def\next{#2}
\ifx\next\empty 
     \setbox\TheAdd=\hbox to\hsize{\hfill\copy0\hfill}
\else\setbox1=\vtop{{\obeylines\hsize=\xs#2\par}}\def\next{#3}
\ifx\next\empty 
     \setbox\TheAdd=\hbox to\hsize{\hfill\copy0\hfill\copy1\hfill}
\else\setbox2=\vtop{{\obeylines\hsize=\xs#3\par}}\def\next{#4}
\ifx\next\empty\ 
     \setbox\TheAdd=\vtop{\hbox to\hsize{\hfill\copy0\hfill\copy1\hfill}
                \vskip20pt\hbox to\hsize{\hfill\copy2\hfill}}
\else\setbox3=\vtop{{\obeylines\hsize=\xs#4\par}}
     \setbox\TheAdd=\vtop{\hbox to\hsize{\hfill\copy0\hfill\copy1\hfill}
	        \vskip20pt\hbox to\hsize{\hfill\copy2\hfill\copy3\hfill}}
\fi\fi\fi\catcode'015=5}}\gdef\Address{\obeylines\GetAddress}

\begin{document}

\Title
The Sphere Packing Problem
\ShortTitle
\SubTitle
\Author
N. J. A. Sloane
\ShortAuthor
\EndTitle
\Abstract
A brief report on recent work on the sphere-packing problem.
\EndAbstract
\MSC
52C17
\EndMSC
\KEY
Sphere packings; lattices; quadratic forms; geometry of numbers
\EndKEY
\Address
N. J. A. Sloane
AT\&T Labs-Research
180 Park Avenue 
Florham Park NJ 07932-0971 USA
njas@research.att.com
\Address
\Address
\Address
\EndAddress
\section{Introduction}
The sphere packing problem has its roots in geometry and number theory
(it is part of Hilbert's 18th problem),
but is also a fundamental question in information theory.
The connection is via the sampling theorem.
As Shannon observes in his classic 1948 paper
\cite{Shan48} (which ushered in the age of digital communication),
if $f$ is a signal of bandwidth $W$ hertz,
with almost all its energy concentrated in an interval of $T$ secs,
then $f$ is accurately represented by a vector of $2WT$ samples,
which may be regarded as the coordinates of a single point in $\RR^n$, $n=2WT$.
Nearly equal signals are represented by neighboring points,
so to keep the signals distinct, Shannon represents them by $n$-dimensional
`billiard balls', and is therefore led to ask:
what is the best way to pack `billiard balls' in $n$ dimensions?

This talk will report on a few selected developments that have taken place
since the appearance of Rogers' 1964 book on the subject, proceeding
upwards in dimension from 2 to 128.
The reader is referred to \cite{SPLAG} (especially the third edition,
which has 800 references covering 1988-1998) for further information,
definitions and references.
See also the lattice data-base \cite{NeSl}.

\section{Dimension 2}
The best packing in dimension 2 is the familiar `hexagonal lattice'
packing of circles, each touching six others.
The centers are the points of the root lattice $A_2$.
The {\em density} $\Delta$ of this packing is the fraction of the plane occupied by the spheres:
$\pi / \sqrt{12} = 0.9069 \ldots$.

In general we wish to find $\Delta_n$,
the highest possible density of a packing of equal nonoverlapping spheres in $\RR^n$,
or $\Delta_n^{(L)}$, the highest density of any packing in which the centers form a lattice.
It is known (Fejes T\'{o}th, 1940) that $\Delta_2 = \DtL = \pi / \sqrt{12}$.
An $n$-dimensional lattice $\Lambda$ of determinant $d$ and minimal nonzero squared
length (or {\em norm}) $\mu$ has packing radius $\rho = \sqrt{\mu} /2$ and
density $\Dt = V_n \rho^n / \sqrt{\det \La}$, where
$V_n = \pi^{n/2} / (n/2)!$ is the volume of a unit sphere.
The {\em center density} of a packing is $\delta = \Dt / V_n$.

We are also interested in
packing points on a sphere, and especially in the `kissing number problem':
find $\tau_n$ (resp. $\tau_n^{(L)}$),
the maximal number of spheres that can touch
an equal sphere in $\RR^n$ (resp. in any lattice in $\RR^n$).
It is trivial that $\tau_2 =  \tau_2^{(L)} = 6$.

\section{Dimension 3}
In spite of much recent work
(\cite{Hal97}, \cite{Hsi93})
$\Dt_3$ is still unknown; nor is $\Dt_n$ known in any dimension above 2.
It is conjectured that $\Dt_3 = \pi / \sqrt{18} =0.74048 \ldots$, as in
the face-centered cubic (f.c.c.) lattice $A_3$.
Muder \cite{Mude93} has shown that
$\Dt_3 \le 0.773055 \ldots$.
It is worth mentioning, however, that there are packings of congruent ellipsoids with 
density considerably greater than $\pi / \sqrt{18}$
\cite{BezKu91}.

In two dimensions the hexagonal lattice 
is (a) the densest lattice packing, (b)~the least dense lattice covering,
and (c)~is geometrically similar to its dual lattice.
There is a little-known three-dimensional lattice that is similar to its dual,
and, among all lattices with this property,
is both the densest packing and the least dense covering.
This is the m.c.c. (or {\em mean-centered cuboidal}) lattice
\cite{CoSl94} with Gram matrix
$$
\frac{1}{2}
\left[
\matrix{
1+ \sqrt{2} & 1 & 1 \cr
1 & 1+ \sqrt{2} & 1- \sqrt{2} \cr
1 & 1- \sqrt{2} & 1 + \sqrt{2} \cr
}
\right] ~.
$$
In a sense this lattice is the geometric mean of
the f.c.c. lattice and its dual the body-centered cubic (b.c.c.) lattice.
Consider the lattice generated by the vectors
$(\pm u, \pm v, 0)$ and $(0, \pm u, \pm v)$ for real
numbers $u$ and $v$.
If the ratio $u/v$ is respectively $1$, $2^{1/2}$ or $2^{1/4}$ we
obtain the f.c.c., b.c.c. and m.c.c. lattices.
The m.c.c. lattice also recently arose in a different context, as the
lattice corresponding to the period matrix of the hyperelliptic Riemann
surface $w^2 = z^8 -1$

\section{Dimensions 4--8}
Table~\ref{ta1} summarizes what is presently known about the sphere packing
and kissing number problems in dimensions $\le$ 24.
Entries enclosed inside a solid line are known to be optimal,
those inside a dashed line optimal among lattices.
\begin{table}[htb]
\begin{center}
\input tab-num.pstex_t
\end{center}
\caption{Densest packings and highest kissing numbers known in low
dimensions.
(Parenthesized entries are nonlattice arrangements that are better than any known lattice.)}
\label{ta1}
\end{table}

The large box in the `density' column refers to Blichfeldt's 1935 result
that the root lattices $\ZZ \simeq A_1$, $A_2$, $A_3 \simeq D_3$,
$D_4$, $D_5$, $E_6$, $E_7$, $E_8$ achieve $\Dt_n^{(L)}$ for $n \le 8$.
It is remarkable that more than 60 years later $\Dt_9^{(L)}$ is still
unknown.

The large box in the right-hand column refers to Watson's 1963
result that the kissing numbers of the above lattices, together
with that of the laminated lattice $\La_9$, achieve
$\tau_n^{(L)}$ for $n \le 9$.
Odlyzko and I \cite[Ch.~13]{SPLAG} and
independently Levenshtein determined $\tau_8$ and
$\tau_{24}$. The packings achieving these two bounds are unique
\cite[Ch.~14]{SPLAG}.

\paragraph{The `Low Dimensional Lattices' project}
Some years ago Conway and I noticed that there were several
places in the literature where the results
could be simplified if they were described in terms of lattices
rather than quadratic forms.
(It seems clearer to say `the lattice $E_8$' rather than `the quadratic form
$2x_1^2 + 2x_2^2 + 4x_3^2 + 4x_4^2 + 20x_5^2 + 12x_6^2 + 4x_7^2 + 2x_8^2 +
2x_1x_2 + 2x_2x_3 + 6x_3x_4 + 10x_4x_5 + 6x_5x_6 + 2x_6x_7 + 2x_7x_8$'.)
This led to a series of papers \cite{CSLDL1}, \cite{CoSl91a}, \cite{CoSl95a}.

Integral lattices of determinant $d=1$ (`unimodular' lattices) have been
classified in dimensions $\le 25$, dimensions 24, 25 being due to
Borcherds. In \cite[Ch.~15]{SPLAG} and \cite[(I)]{CSLDL1}
we extended this to $d \le 25$ for various ranges of dimension.

\cite[(II)]{CSLDL1} is based on the work of Dade,
Plesken, Pohst and others,
and describes the lattices associated with the maximal irreducible subgroups of 
$GL (n, \ZZ )$ for $n=1, \ldots, 9, 11, 13, 17, 19, 23$.
Nebe, and Nebe and Plesken (see \cite{Nebe96a}, \cite{Plesk96})
have recently completed the enumeration of the maximal finite
irreducible subgroups of $GL(n, \QQ )$ for $n \le 31$,
together with the associated lattices.

\cite[(IV)]{CSLDL1} gives an improved version of
the mass formula for lattices, and
\cite[(V)]{CSLDL1} studies when an $n$-dimensional
integral lattice can be represented as 
a sublattice of $\ZZ^m$ for some $m \ge n$,
or failing that, by a sublattice of $s^{-1/2} \ZZ^m$ for some integer $s$.
\cite{CoSl91a} describes the Voronoi and Delaunay cells of all the root lattices
and their duals, and
\cite[(VI),~(VIII)]{CSLDL1} discusses how the Voronoi cell of a 3- or 4-dimensional lattice
changes as the lattice is continuously varied.

\cite[(VII)]{CSLDL1} determines the `coordination sequences' of various lattices.
Consider $E_8$, for example, and let $S(k)$ denote the number of lattice points
that are $k$ steps from the origin, where a step is a move to an adjacent sphere
($S(1)$ is the kissing number).
Then $\sum_{k=0}^\infty S(k) x^k = f(x) / (1-x)^8$, where
$f(x) = 1+ 232 x + 7228x^2 + \ldots + x^8$.
Thus the coordination sequence for $E_8$ begins 1, 240, 9120, $\ldots$.
For other examples see \cite{SloEIS}

\paragraph{Perfect lattices}
One possible approach to the determination
of the densest lattices in dimensions 7 to 9
is via Voronoi's theorem that the density of
$\La$ is a local maximum if and only if $\La$ is perfect and eutactic \cite{Mar96}.

In 1975 Stacey, extending the work of several earlier
authors, published a list of 33 perfect lattices in dimension 7.
Unfortunately one of the 33 was omitted from her papers and her dissertation.
In \cite[(III)]{CSLDL1} we reconstructed the missing lattice and `beautified'
all 33, computing their automorphism groups, etc.
In 1991 Jaquet-Chiffelle \cite{Jaq93} completed this
work by showing that this is indeed the full list
of perfect lattices in $\RR^7$.
This provides another proof that $E_7$ is the densest lattice in dimension 7.

Martinet, Berg\'{e} and their students are presently attempting to classify
the eight-dimensional perfect lattices, and it appears that there will be
roughly 10000 of them.
Whether this approach can be used to determine $\Dt_9^{(L)}$ remains to be seen!

\section{Dimension 9. Laminated lattices}
\begin{figure}[htb]
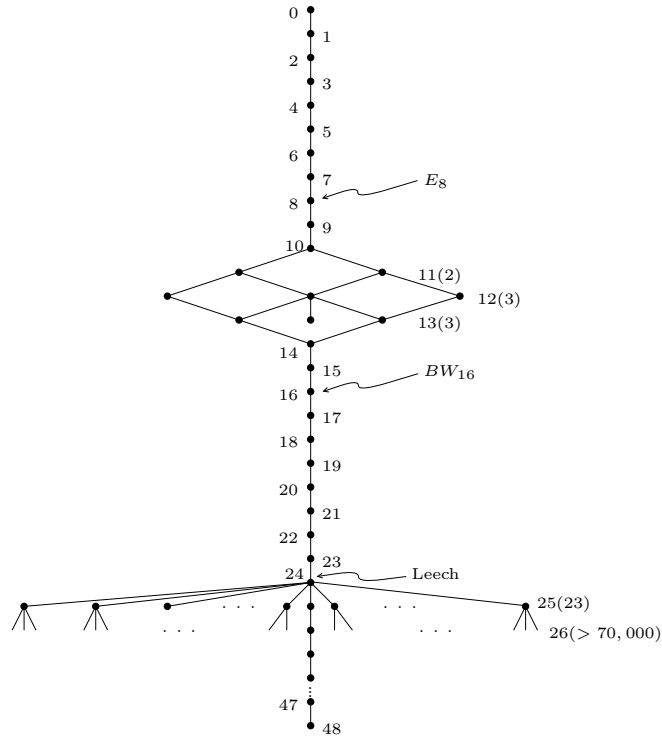

\begin{center}
\input dots.pstex_t
\end{center}
\caption{Inclusions among laminated lattices $\La_n$.}
\label{fg2}
\end{figure}

There is a simple construction, the `laminating' or `greedy' construction,
that produces many of the densest lattices in dimensions up to 26.
Let $\La_1$ denote the even integers in $\RR^1$, and define the $n$-dimensional laminated lattices
$\La_n$ recursively by:
consider all lattices of minimal norm 4 that contain some $\La_{n-1}$ as a sublattice, and select those of greatest density.
It had been known since the 1940's that this produces the densest lattices 
known for $n \le 10$.
In \cite{Con32} we determined {\em all}
inequivalent laminated lattices for $n \le 25$, and found the density of $\La_n$ for
$n \le 48$ (Fig. \ref{fg2}).
A key result needed for this was the determination of the covering radius
of the Leech lattice and the enumeration of the deep holes in that lattice
\cite[Ch.~23]{SPLAG}.

\paragraph{What are all the best sphere packings in low dimensions?}
In \cite{CoSl95a} we describe what may be {\em all} the best
packings in dimensions $n \le 10$,
where `best' means both having the highest density and not permitting
any local improvement.
In particular, we conjecture that $\Dt_n^{(L)} = \Dt_n$ for $n \le 9$.
For example, it appears that
the best five-dimensional sphere packings are parameterized by
the 4-colorings of $\ZZ$.
We also find what we believe to be the exact numbers of 
`uniform' packings among these,
those in which the automorphism group acts transitively.
These assertions depend on certain plausible but as yet unproved postulates.

\paragraph{A remarkable property of 9-dimensional packings.}
We also show in \cite{CoSl95a} that the laminated lattice $\La_9$
has the following astonishing property.
Half the spheres can be moved bodily through arbitrarily large distances without overlapping the other half,
only touching them at isolated instants,
the density remaining the same at every instant.
A typical packing in this family consists of the points of
$D_9^{\theta +} = D_9 \cup D_9 + ((1/2)^8, \theta /2)$,
for $\theta$ real.
$D_9^{0+}$ is $\La_9$ and $D_9^{1+}$ is $D_9^+$, the
9-dimensional diamond structure.
All these packings have the same density, which we conjecture is the value of $\Dt_9 = \Dt_9^{(L)}$.
Another result in \cite{CoSl95a} is that
there are extraordinarily many 16-dimensional packings that are just as dense
as the Barnes-Wall lattice $BW_{16} \simeq \Lambda_{16}$.

\section{Dimension 10. Construction A.}
In dimension 10 we encounter for the first time a nonlattice packing that is denser than all known lattices.
This packing, and the nonlattice packing with the highest known kissing
number in dimension 9,
are easily obtained from `Construction A'
(cf. \cite{Lee10}).
If ${\cal C}$ is a binary code of length $n$, the corresponding packing is
$P( {\cal C} ) = \{ x \in \ZZ^n: x ~(\bmod~2) \in {\cal C} \}$.

Consider the vectors $abcde \in ( \ZZ / 4 \ZZ)^5$ where $b,c,d \in \{ +1, -1\}$,
$a= c-d$, $e= b+c$, together with all their cyclic shifts, and apply the
`Gray map'
$0 \to 00 , ~ 1 \to 01, ~ 2 \to 11 , ~ 3 \to 10$
to obtain a binary code ${\cal C}_{10}$ containing 40 vectors of length 10 and
minimal distance 4.
This is our description \cite{CoSl94b} of a code first discovered by Best.
The code is unique \cite{LiV94}.
Then $P( {\cal C}_{10} ) = P_{10c}$ is the record 10-dimensional packing.

\begin{figure}
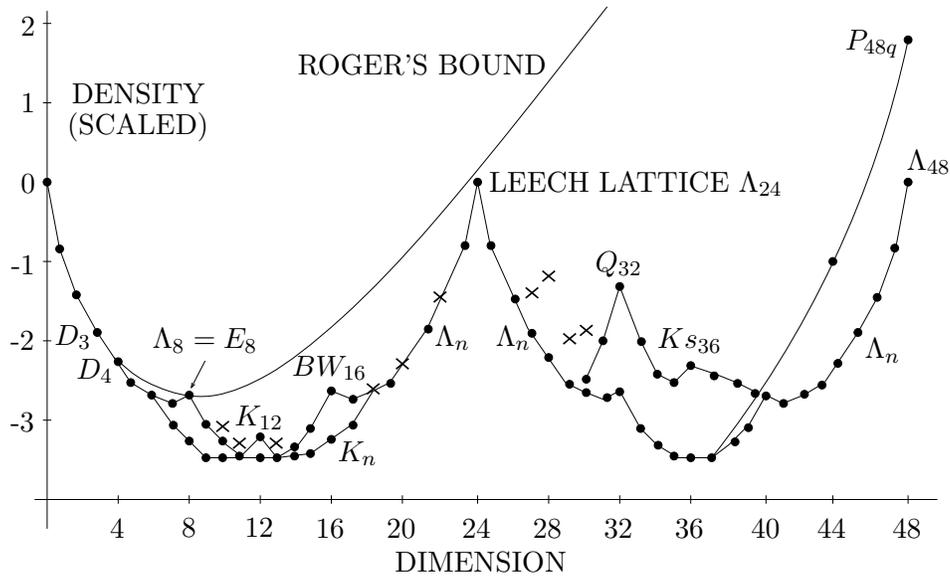

\begin{center}
\input curve.pstex_t
\end{center}
\caption{Densest sphere packings known in dimensions $n \le 48$.}
\label{fg3}
\end{figure}

Figure \ref{fg3} shows the density of the best packings known up to dimension 48,
rescaled to make them easier to read.
The vertical axis gives $\log_2 \delta + n(24-n)/96$.
The figure also shows the upper bounds of Muder (for $n=3$) and Rogers $(n \ge 4)$.
Lattice packings are indicated by small circles, nonlattices by crosses
(however, the locations of the lattices are only approximate).
The figure is dominated by the two arcs of the graph of the laminated lattices $\La_n$,
which touch the zero ordinate at $n=0$, 24 (the Leech lattice) and 48.
$K_{12}$ is the Coxeter-Todd lattice.

\section{Dimensions 18--22}
Record nonlattice packings in dimensions 18, 20 and 22 have recently
been given in \cite{BiEd98}, \cite{CoSl96}, \cite{Vard95}.
Vardy's construction \cite{Vard95},
`Construction $B^\ast$',
also uses binary codes.
Let ${\cal B}$ and ${\cal C}$ be codes of length $n$ such that
$c \cdot ( {\bf 1} + b) = 0$ for all $b \in {\cal B}$, $c \in {\cal C}$,
and set
$P^\ast ( {\cal B} , {\cal C} ) = \{ {\bf 0} + 2b + 4x , {\bf 1} + 2c+ 4y ~:~
b \in {\cal B} , c \in {\cal C} , x,y \in \ZZ^n , \sum x_i ~{\rm even} ,
\sum y_i ~{\rm odd} \}$.
For example, by taking ${\cal B}$ to be the quadratic residue code of length 18
and ${\cal C}$ to be its dual,
Bierbrauer and Edel \cite{BiEd98} obtain a new record packing in $\RR^{18}$.

\section{Dimension 24. The Leech lattice}
The Leech lattice $\La_{24}$ is a remarkably dense packing
in $\RR^{24}$ (as can be seen from Fig.~\ref{fg3}).
Here are four constructions.
(i)~As a laminated lattice: start in dimension 1 with the lattice
$\La_1 = \ZZ$ and apply the greedy algorithm (see Fig.~\ref{fg2}).
(ii) Apply Construction A to the Golay code of length 24 to obtain a lattice
$L_{24}$.
Then $\La_{24}$ is spanned by $(-3/2, 1/2, \ldots, 1/2 )$ and $\{ x \in L_{24} : \sum x_i \equiv 0$ $(\bmod~4)\}$.
(iii)~Hensel lift the Golay code to an extended cyclic (and self-dual)
code over $\ZZ / 4 \ZZ$ and apply
`Construction A mod 4' \cite{BonCS95}.
(iv)~There is a unique unimodular even lattice ${\rm II}_{25,1}$
in Lorentzian space $\RR^{25,1}$, consisting of the points
$(x_0 x_1 \cdots x_{24} | x_{25} )$ with all $x_i \in \ZZ$ or all
$x_i \in \ZZ + 1/2$ and satisfying $x_0 + \cdots + x_{24} - x_{25} \in 2 \ZZ$.
Let $w = (0~ 1 \cdots 24 | 70)$,
a vector of zero length.
Then $(w^\perp$ in ${\rm II}_{25,1} ) / w$ is $\La_{24}$
\cite[Ch.~26]{SPLAG}.

\section{Dimensions 26--31}
New packings in these dimensions have been discovered by Bacher, Borcherds,
Conway, Vardy, Venkov --- see \cite{SPLAG} for details.

\section{Dimension 32. Modular lattices}
An {\em $N$-modular} lattice \cite{Queb95} is an integral lattice that is similar to its dual,
under a similarity that multiplies norms by $N$.
A unimodular lattice is 1-modular.
The interest in this family arises because many of the densest
known lattices are $N$-modular:
$\ZZ$, $A_2$, $D_4$, $E_8$, $K_{12}$, $BW_{16}$, $\La_{24}$, $Q_{32}$, $P_{48q}$, $\ldots$.

Quebbemann's lattice $Q_{32}$, for example, is 2-modular,
and can be constructed from a Reed-Solomon code of length 8 over $\FF_9$
\cite{Que6}, \cite[Ch.~8]{SPLAG}.

\paragraph{Shadow theory.}
The concept of the shadow of a lattice or code was introduced
in \cite{CoSl90}, \cite{CoSl90a} (see also \cite{CoSl98}) and has proved
to be very useful (\cite{CoSl90a} has stimulated over 50 sequels in the
coding literature).

Let $\La$ be an $n$-dimensional unimodular lattice.
If $\La$ is even then the
{\em shadow} $S(\La ) = \La$, otherwise
$S( \La ) = ( \La_0)^\ast \setminus \La$, where the subscript 0 denotes even sublattice.
The set $2S( \La ) = \{ 2s : s \in S( \La ) \}$ is precisely the set of {\em parity vectors} for $\La$,
i.e. the vectors $u \in \La$ such that $u \cdot x \equiv x \cdot x$ $(\bmod 2)$ for all $x \in \La$.
Such vectors have been studied by many authors
from Braun (1940) onwards,
but their application to obtaining bounds on lattices seems to have been overlooked.

If the theta series of $\La$ is $\Theta_\La (z)$ then \cite{CoSl90} the shadow has theta series
\begin{equation}\label{Eq1}
\left( \frac{e^{\pi i/4}}{\sqrt{z}} \right)^n \Theta_\La \left( 1- \frac{1}{z} \right) ~.
\end{equation}
One of the most satisfying properties of integral lattices is the classical theorem that 
(a) if $\La$ is a
unimodular lattice then $\Theta_\La$ belongs to the graded ring
$\CC [ \Theta_{\ZZ} , \Theta_{E_8} ]$,
and (b) if $\La$ is even then $\Theta$ belongs to
$\CC [ \Theta_{E_8} , \Theta_{\La_{24}} ]$.

To illustrate the use of the shadow, let us prove there is no 9-dimensional
unimodular lattice of minimal norm 2.
If so then from (a)
$\Theta_\La = - \Theta_{\ZZ}/8 + 9 \Theta_{E_8}/8 = 1+ 252 q^2 + 456 q^3 + \cdots$,
where $q = e^{\pi iz}$.
But then (\ref{Eq1})
implies $\Theta_{S( \La )} = \frac{9}{4} q^{1/4} + \frac{1913}{4} q^{9/4} + \cdots$,
a contradiction since $\Theta_{S(\La )}$ must have integer coefficients.

In \cite{Mal1} we used (a), (b) to show that the minimal norm
$\mu$ of an $n$-dimensional odd unimodular lattice satisfies
\begin{equation}\label{Eq4}
\mu \le \left[ \frac{n}{8} \right] + 1 ~,
\end{equation}
and for an even unimodular lattice
\begin{equation}\label{Eq5}
\mu \le 2 \left[ \frac{n}{24} \right] +2 ~.
\end{equation}
In \cite{RaSl98a}
we used shadow theory to strengthen (\ref{Eq4}) by showing that odd lattices satisfy
\begin{equation}\label{Eq6}
\mu \le 2 \left[ \frac{n}{24} \right] +2 ~,
\end{equation}
except that $\mu \le 3$ when $n =23$.
In view of the similarity between (\ref{Eq5}) and (\ref{Eq6}) we propose that a lattice
satisfying either bound with equality be called {\em extremal}
(the old definition of this term was based on (\ref{Eq4}) and (\ref{Eq5})).

Quebbemann \cite{Queb97} has generalized (\ref{Eq5}) to certain families of even $N$-modular lattices, and
analogous bounds for odd $N$-modular lattices (using an appropriate generalization of the 
shadow) were given in \cite{RaSl98a}.
One can then define extremal $N$-modular lattices.

\section{Higher dimensions}
Space does not permit more than a mention of the following:
Kschischang and Pasupathy's lattice $Ks_{36}$ in $\RR^{36}$ \cite{KsP92};
the three extremal unimodular lattices $P_{48q}$, $P_{48p}$,
$P_{48n}$ in $\RR^{48}$,
the latter being a recent discovery of Nebe \cite{Nebe98};
Bachoc's extremal 2-modular lattice in $\RR^{48}$ \cite{Baco97};
Nebe's extremal 3-modular lattice in $\RR^{64}$ \cite{Nebe98};
and Bachoc and Nebe's extremal unimodular lattice in $\RR^{80}$ \cite{BacoN98}.

The existence of the following extremal lattices is an open question:
3-modular in $\RR^{36}$ (determinant $d= 3^{18}$, minimal norm $\mu =8$);
2-modular in $\RR^{64}$ $(d= 2^{32}, \mu =10 )$;
unimodular in $\RR^{72}$ ($d=1$, $\mu =8$).

From dimensions 80 to about 4096 the densest lattices known are the Mordell-Weil 
lattices discovered by Elkies \cite{Elki94},
and Shioda \cite{Shiod91}.
But we know very little about this range, as evidenced by the recent 
construction
of record kissing numbers in dimensions 32 to 128 \cite{EdRS98} 
from binary codes.
In dimension 128, for example, the Mordell-Weil lattice has kissing
number 218044170240 \cite{Elki},
whereas in our construction
(which admittedly is not a lattice)
some spheres touch 8812505372416 others.

It would also be desirable to have better upper bounds,
especially in low dimensions (see Fig.~\ref{fg3}).
The Kabatiansky-Levenshtein bound
is asymptotically better than the Rogers' bound, but not
until the dimension is above about 40.
We know very little about these problems!

In short, many beautiful packings have been discovered, but there are few proofs that any of them are optimal.

\Addresses
\end{document}